\documentclass[12pt]{article}
\usepackage{amssymb,amsmath,amsthm,bbm,verbatim}
\newtheorem{corollary}{Corollary}[section]
\newtheorem{lemma}[corollary]{Lemma}
\newtheorem{proposition}[corollary]{Proposition}

\newtheorem{theorem}[corollary]{Theorem}
\theoremstyle{definition}

\newtheorem*{remark}{Remark}

\newcommand{\Prob} {{\mathbb P}}
\newcommand{\Z}{{\mathbb Z}} 
\newcommand{\E}{{\mathbb E}}

\newcommand{\R}{{\mathbb{R}}}
\newcommand{\C}{{\mathbb C}}

\newcommand{\dist}{{\rm dist}}

\def \p {\partial}

\def \loops{{\cal L}}

\def \walkmeasure {{\nu}}

\def  \brownmeas  {{\mu^{\rm u}}}
\def \bridgemeas {{\mu^{\rm br}}}
\def \rootmeas {{\mu}}
\def \Disk {{\mathbb D}}
\def \walkmeasureroot {{\mu^{\rm rw}}}

\begin{document}

\vspace{1.5in}

\begin{center}
\begin{LARGE}
{\bf Random Walk Loop Soup}\\
\vspace{2ex}
\end{LARGE}
\begin{large}
Gregory F. Lawler\footnote{Research supported by the National Science
Foundation}\\
Jos\'e A. Trujillo Ferreras\\
\vspace{6ex}
\end{large}
Department of Mathematics\\
Malott Hall\\
 Cornell University\\
Ithaca, NY 14853-4201\\
lawler@math.cornell.edu\\
jatf@math.cornell.edu
 \end{center}

\vspace{.5in}

\begin{abstract}
The Brownian loop soup introduced in \cite{LW} is a Poissonian
realization from a $\sigma$-finite measure on unrooted
loops.  This measure satisfies both conformal invariance
and a restriction property.  In this paper, we define
a random walk loop soup and show that it converges to
the Brownian loop soup.  In fact, we give a strong approximation
result making use of the strong approximation result
of Koml\'os, Major, and Tusn\'ady.  To make the
paper self-contained, we include a proof of the approximation
result that we need.
\end{abstract}

\section{Introduction}

The Brownian loop soup with intensity $\lambda$, which we define
below,  is a Poissonian
realization from a particular measure on unrooted
loops in $\C$  that satisfies
both conformal invariance and a property called the restriction property. 
A realization of the loop soup consists of a countable collection
of loops. In a fixed bounded domain $D$, there are an infinite
number of loops that stay in $D$; however, the number of loops
of diameter at least $\epsilon$ in the bounded domain is
finite.  A corollary of conformal invariance is scale invariance:
 if ${\cal A}$ is a realization of the Brownian
loop soup and each loop is scaled in space by $1/N$ and in time
by $1/N^2$, the resulting configuration also has the distribution
of the Brownian loop soup. 
In this paper, we will show that the Brownian soup is a limit of 
random walk soups.  There are two natural approaches to showing this.
One is a ``weak'' limit to show that the Brownian loop measure is
a weak limit of random walk measures (this requires some care since
the measures are infinite).  However, we choose
the more direct ``coupling'' approach
of defining the random walk loop soup and the Brownian
loop soup on the same probability space
so that the realizations are close.  Since a realization
is a countable collection of loops, it is a little tricky to say what
it means for the realizations to be close.  We will prove, in fact,
that in a bounded domain $D$, except for an event of small probability,
there is a one-to-one correspondence
between the Brownian loops and the random walk loops if we restrict
to loops that are not ``too small''.
The Brownian loops and random walk loops that correspond with each
other will be very close.
 We will use the dyadic
 approximation scheme as in \cite{KMT}
 to establish the strong approximation
of the two soups.

We start by defining the Brownian loop measure.
It is easier to define the loop measure first on {\em rooted}
loops.  A {\em (rooted) loop} is a continuous function
$\gamma:[0,t_\gamma] \rightarrow \C$ with $\gamma(0) =
\gamma(t_\gamma)$.  We will only consider loops with
$0 < t_\gamma < \infty$.  Let ${\cal C}$ denote the set of
all loops and ${\cal C}_t$ the set of loops $\gamma$ with
$t_\gamma = t$ and $\gamma (0) = \gamma (t_\gamma)= 0$.  The {\em Brownian bridge} measure
$\bridgemeas$ is the probability measure on loops induced
by a Brownian bridge , i.e., by $B_t := W_t -t W_1, 0 \leq t \leq 1,$
where $W_t$ is a standard two-dimensional Brownian motion.
The measure $\bridgemeas$ is supported on ${\cal C}_1$.  
The {\em (rooted) Brownian loop measure} is the measure $\rootmeas$
on ${\cal C}_1 \times \C \times (0,\infty)$ given by
\[   
\rootmeas= \bridgemeas \times {\rm area} \times [\frac{1}{2\pi t^2} \, dt] . \]
The measure $\rootmeas$    induces a measure on ${\cal C}$, which
we also denote by $\rootmeas$, by the function $(\gamma,z,t) \mapsto \tilde \gamma$,
where $\tilde \gamma$ is $\gamma$ scaled (using Brownian scaling) to
have time duration $t$ and translated to have root $z$.  In other
words,
\[    \tilde \gamma(s) = z + t^{1/2} \, \gamma(s/t), \;\;\;\;
   0 \leq s \leq t . \]
This measure is clearly translation invariant, and it is straightforward
  to check that if $ r > 0$, then $\rootmeas$ is invariant
under the Brownian scaling
map  $(\gamma,z,t) \mapsto (\gamma,rz,r^2t)$.

Let us denote by $\mu^{\rm br}_t(z)$ the probability measure on loops induced 
by a Brownian bridge of time duration $t$ rooted at $z$. Then the measure 
$\rootmeas$ (as a measure on ${\cal C}$) can be written as
\[ \int _\C \int _0 ^\infty \frac{1}{2\pi t ^2} \, \mu_t^{{\rm br}}(z) \, dt \,dz.
\]

An {\em unrooted loop} is an equivalence class of (rooted) loops under the equivalence
$\gamma \sim \theta_r \gamma$ for every $r  \in \R$, where
$\theta_r \gamma(s) = \gamma(s+r)$ (here we consider a 
rooted  loop $\gamma$ of time duration $t_\gamma$
 as a continuous function $\gamma: \R
\rightarrow \C$ with $\gamma(s+t_\gamma) = \gamma(s)$ for all $s$).
The unrooted loop measure $\brownmeas$ is the measure obtained from $\rootmeas$
by ``forgetting the root.''  A {\em rooted Brownian loop soup}
with intensity $\lambda$ is a Poissonian realization from $\lambda \rootmeas$.
An {\em (unrooted) Brownian loop soup} is a realization from $\lambda
\brownmeas$.  One can obtain an unrooted loop soup by starting with
a rooted loop soup and forgetting the root.

If $D$ is a domain in $\C$ we let $\rootmeas_D,\brownmeas_D$ denote $\rootmeas,
\brownmeas$ restricted to loops that lie in $D$.  The
family of measure $\{\brownmeas_D\}$ clearly satisfy the restriction
property, i.e., if $D' \subset D$ then $\brownmeas_{D'}$ is $\brownmeas_D$
restricted to curves lying in $D'$.  It is also shown in \cite{LW} that
the family satisfies a conformal invariance property, i.e., if
$f:D \rightarrow D'$ is a conformal transformation, then
$f \circ \brownmeas_D = \brownmeas_{D'}$, if the quantites are suitably
interpreted.  In particular, if $\gamma$ is a curve lying in $D$, we
define $f \circ \gamma$ to be the curve in $D'$,  reparametrized
by the conformal map; see \cite{LW} for details.  The measure $\mu$
on rooted loops is not conformally invariant.

In this paper we study the loop measure for simple (nearest
neighbor) random walks on the integer lattice $\Z^2$,  which we can
consider as a subset of $\C$.
The rooted loop measure $\walkmeasureroot$ gives each 
(nearest neighbor) random walk loop in $\Z^2$ of length $2n$
 measure $(2n)^{-1} \, 4^{-2n}$.   The unrooted loop measure
is obtained from the rooted loop measure by ``forgetting the root''.
It is almost the same measure as that obtained by giving measure $4^{-2n}$
to every unrooted loop
of length $2n$. (If a loop of length $2n$ is obtained by taking
a loop of length $n$ and repeating the same loop again, then this unrooted
loop
does not get full measure $4^{-2n}$ under our random walk loop
measure; these exceptional loops are an
exponentially small subset of the set of all loops so it is not important
whether we give these unrooted loops measure $4^{-2n}$ or $(1/2) \, 4^{-2n}$.)
We will focus on the rooted measure in this paper.
 A rooted random walk loop of length $2n$ can also be considered as
a continuous path $\gamma:[0,2n] \rightarrow \C$ by linear
interpolation.
 We will call
a Poissonian realization from $\lambda\, \walkmeasureroot$ a
{\em rooted random walk loop soup} (with intensity
$\lambda$).

In this paper, we make a precise statement that the random walk loop
soup, appropriately scaled, approaches the Brownian loop soup. 
We will define $({\cal A}_\lambda,
\tilde {\cal A}_\lambda)$ on the same probability space  
 so that ${\cal A}_\lambda$ is a realization
of the Brownian loop soup with intensity $\lambda$,
 and $\tilde {\cal A}_\lambda$
is a realization of the random walk loop soup with intensity
$\lambda$. We consider the loops in $\tilde {\cal A}_\lambda$
as curves in ${\cal C}$ by linear interpolation.
  Note that ${\cal A}_\lambda$ is a (random) countable set
of curves and $\tilde {\cal A}_\lambda$ is a (random) multi-set
(i.e.,
a set where some elements can appear more than once) of
lattice curves.  For each
positive integer $N$ we define ${\cal A}_{\lambda,N}$ to be
the collection of loops obtained from ${\cal A}_\lambda$ by
scaling space by $1/N$.  More precisely,
\[  {\cal A}_{\lambda,N} = \{ \Phi_N\gamma: \gamma \in {\cal A}_\lambda \} , \]
where $t_{\Phi_N \gamma} =  t_\gamma/N^2$ and
\[          \Phi_N \gamma(t) = N^{-1} \gamma(tN^2), \;\;\;\;
     0 \leq t \leq t_\gamma/N^2. \]
Note that the scaling rule implies that ${\cal A}_{\lambda,N}$ is
a realization of the Brownian loop soup with parameter $\lambda.$
We define
\[   \tilde A_{\lambda,N} = \{\tilde \Phi_N \gamma: \gamma
     \in \tilde {\cal A}_\lambda\} , \]
where $t_{\tilde {\Phi}_N \gamma} = t_{\gamma}/(2N^2)$ and
\[    \tilde \Phi_N \gamma(t) = N^{-1} \gamma(t 2N^2),
\;\;\;  0 \leq t \leq t_{\gamma}/(2N^2). \] 
  The
scaling is slightly different for the random walk loops because
the covariance of a simple two-dimensional random walk in
$2n$ steps is $nI$ as opposed to $2nI$ for a Brownian motion
at time $2n$; roughly speaking, this is because in $2n$ steps,
the random walk moves about $n$ steps horizontally and 
$n$ steps vertically.

 We will prove the theorem below. The ideas in the proof are simple 
and flexible. However, due to discretization, stating the result is 
a little unwieldy. To aid, we introduce the following auxiliary functions. 
For $t \ge (5/8) N ^ {-2}$ and positive integer $k$,  we let 
\[  \varphi _N ( t ) = \frac{2k}{2 N ^2} \;\;\;\;
\mbox{if} \;\;\;\;
 \frac k {N ^2} - \frac 3 {8 N ^2} \le t < \frac k {N^2} + \frac 5 {8 N ^2}.\]
Also, for $z \in \C, z_0 \in \Z^2$ we define 
\[ \psi _N ( z ) = \frac {z _0}N\;\;\;\; \mbox{if}
\;\;\;\;
\max \{|\rm{Re} \{ Nz -z_0 \}|,|\rm{Im} \{ Nz - z _0 \}| \} < \frac 1 2.\]
The definition of $\psi _N$ if $Nz$ happens to fall on
a bond of the dual lattice of $\Z ^2$ is irrelevant for our theorem.

\begin{theorem}  \label{maintheorem2}
One can define on the same probability space ${\cal A}_\lambda$
and $\tilde {\cal A}_{\lambda}$ such that:
\begin{itemize}
\item  For each $\lambda > 0$, ${\cal A}_\lambda$ is a realization
of the Brownian loop soup; the realizations are increasing in $\lambda$.
\item  For each $\lambda > 0$, ${\tilde A}_\lambda$ is a realization
of the random walk loop soup; the realizations are increasing
in $\lambda$,
\end{itemize}
and such that the following holds.
Let ${\cal A}_{\lambda,N}$ and $\tilde {\cal A}_{\lambda,N}$
be as defined above.
 Then there exists a $c>0$ such that for every $r\geq 1,
N, \lambda$ and every $2/3 < \theta <2$,
 except perhaps on an event of probability at most
$c \, (\lambda+1) \, r^2 \, N^{2-3\theta} $, 
  there is a one-to-one
correspondence between $\{ \tilde{\gamma} \in \tilde{\cal A} _{\lambda, N} 
:t _ {\tilde{\gamma} } > N ^ {\theta - 2}, | \tilde{\gamma} (0) | < r \}$ and 
$\{ \gamma \in {\cal A} _{\lambda, N }: 
\varphi _N (t _{\gamma}) > N ^{ \theta - 2 },
| \psi _N (\gamma (0)) | < r \}$.
If $\tilde \gamma \in \tilde{\cal A}_{\lambda,N}$ and 
$\gamma \in {\cal A}_{\lambda,N}$ are paired
in this correspondence, then 
\[   |t_\gamma  - t_{\tilde \gamma}| \leq 5/8 N ^{-2} \]
\[   \sup_{0 \leq s \leq 1}
     |\gamma( s t_\gamma) - \tilde \gamma(s t_{\tilde \gamma})|
                \leq c \, N^{-1}\log N . \]
\end{theorem}

The outline of the paper is as follows.  In the next three sections,
we define the random walk loop soup, state the strong approximation
result between random walk bridges and Brownian bridges that we
need, construct the probability space on which both
the random walk and Brownian loop soups are defined, and verify
that the construction satisifes Theorem \ref{maintheorem2}.
The next section concerns the soups in bounded domains.  Here we
establish a similar result to the theorem above, although the
error terms are somewhat larger.  The remainder of the paper gives
a self-contained proof of the strong approximation result that
we need.

\section{Random walk loops and random walk soup} \label{walksoup}

A {\em (rooted) random walk loop} of length $2n$ in $\Z^2$ (which will
be considered as a subset of $\C$) is 
a $(2n+1)$-tuple $\omega = [\omega_0,\ldots,\omega_{2n}]$
with $|\omega_{j} - \omega_{j-1}| = 1$ and $\omega_0 = 
\omega_{2n}$. A loop can be considered
as a curve $\gamma:[0,2n] \rightarrow \C$; here $\gamma(m)
= \omega_m$ for integer $m$ and $\gamma(t)$ is defined for other
$t$ by linear interpolation.  Let $\loops_n$ denote the set of 
random walk loops of
length $2n$ and $\loops_n^z$ be the set of such loops
with $\omega_0 = z$.  There is a natural one-to-one
correspondence between $\loops_n^0$ and $\loops_n^z$
given by $\omega \leftrightarrow z + \omega$.   We let
$\loops = \cup_{n \geq 1} \loops_n$, 
$\loops^z = \cup_{n \geq 1} \loops_n^z$.
 We let $\walkmeasure$ denote the
random walk loop measure  $\loops$, i.e.,
  the measure that assigns
measure $4^{-2n}$ to each $\omega \in \loops_{2n}$.  We write
$\walkmeasure^z_n$ for $\walkmeasure$ restricted to 
$\loops^z_n$.  For
fixed $z$, it is straightfoward to show
\[          \walkmeasure(\loops^z_n) = 
\walkmeasure(\loops^0_n) = \left[2^{-2n}\,\binom{2n}{n}
  \right]^2 =  \frac{1}{\pi \, n}
            - \frac 1 {4\pi n^2}   + O(\frac{1}{n^3}) . \] 
The second equality follows
from   
the fact that the probability that a two-dimensional
simple random walk returns to the origin at time $2n$ is the
square of the probability that a one-dimensional simple random
walk returns to the origin; see (\ref{twofromone}) below
to see why this is true.
The final equality is  derived from Stirling's
forumla with error:
\begin{equation}  \label{stirling}
 n! =\sqrt{2 \pi} \,  n^{n+(1/2)} \, e^{-n}  \, \left[1 +
  \frac{1}{12n} + O(\frac{1}{n^2}) \right] . 
\end{equation}
 
We define the {\em (rooted) random walk loop measure}
$\walkmeasureroot$ to be the measure that assigns
measure $(2n)^{-1} \walkmeasure(\omega)
= (2n)^{-1} \, 4^{-2n}$ to each
$\omega \in \loops_n$. A {\em rooted random walk loop soup of intensity $\lambda$}
is a Poissonian realization from the measure $\lambda \, \walkmeasureroot$.
We will  obtain  a realization by first
using a Poisson point process to select
a multi-set of 
ordered pairs $(n,z)$, where the length of the
loop is $2n$ and it is rooted at $z$.
Then, given $(n,z)$, we
 choose a loop from the appropriate
random walk bridge measure.  In other words, we write
\[ \walkmeasureroot = \sum_{z \in \Z^2} \sum_{n=1}^\infty
     \frac{\walkmeasure_n^z}{2n}  =
\sum _{z \in \Z ^2} \sum_{n =1} ^\infty \frac{ \walkmeasure(\loops^z_n) }
{2n} \; 
 \frac{\nu _n ^z}{\walkmeasure(\loops^z_n)}.
\]

  Let
\[   \tilde N(n,z;t) , \;\;\;\;n \in \{1,2,3,\ldots\}, \; z \in \Z^2 , \]
be independent Poisson processes (in the variable
$t$) with parameter
\[    \tilde q_n := \frac{1}{2n}  \walkmeasure(\loops^z_n) 
                     =  \frac{1}{2 \pi n^2}  - \frac
  1 {8 \pi n^3}+ O(\frac{1}{n^4}) . \]
Let
\[ \tilde   L(n,z;m), \;\;\;\; n \in \{1,2,3,\ldots\}, \; z \in \Z^2, \;
    m \in \{1,2,3,\ldots\}  \]
be independent random variables, independent of the
$\tilde N(n,z;t)$,  taking values in $\loops^0$; the
distribution of $\tilde L(n,z;m)$ is the probability measure of
a random walk bridge of $2n$ steps, i.e.,
${\nu _n ^0}/{\walkmeasure(\loops^0_n)}$, which is the
  uniform probability measure on $\loops_n^0$. 
If $\omega \in \loops_n^z$,
let
\[
J_t(\omega) = \sum_{k=1}^{\tilde N(n,z;t)}
\mathbbm{1}_{\{L(n,z;k) + z=  \omega\}}. 
\]
Note that $\{J_t(\omega): \omega \in \loops\}$
is a collection of independent Poisson processes; the process
$J_t(\omega)$ has parameter $(2n)^{-1}4^{-2n}$ if $\omega \in \loops_n$.
 We could equally well
have constructed the loop soup starting with these Poisson processes.
We have chosen the longer construction because it will be
useful for coupling the loop soup with the Brownian loop
soup.
 
Although we have used $t$ for the time parameter of the Poisson
processes, by choosing $t = \lambda$ we get an increasing
family of realizations of the loop soup $\tilde {\cal A}_\lambda$
parametrized by
$\lambda$.
We think of the loop soup of intensity
$\lambda$  as a multi-set $\tilde {\cal A}_\lambda$
of loops where loop $\omega$ appears $J_\lambda(\omega)$ times
in $\tilde{\cal A}_\lambda$.

\section{Strong Approximation} \label{strongapprox}

If $S_n$ denotes a simple random walk,
we define $S_t,\, 0 \leq t < \infty$,  by linear interpolation.
The key to the coupling is the following result, due to Koml\'os,
Major, and Tusn\'ady, which shows that a simple random walk bridge and a
Brownian bridge can be coupled very closely.  Because the
form of the result we need is slightly different than that
proved in \cite{KMT}, we have included a proof in the final
section. By a simple random walk bridge of time duration 
$2n$ we will mean a process $X_t, \, 0 \le t \le 2n$, that 
has the law of $S_t,\, 0 \le t \le 2n$ conditioned to have 
$S_{2n} = 0$.

\begin{lemma} [Dyadic approximation] 
There exists a $c < \infty$ such that  for every positive integer
$n$, there exists a probability space $(\Omega,{\cal F},\Prob)$
on which are defined a one-dimensional
Brownian bridge $B_t, 0 \leq t \leq 1$
and a one-dimensional simple
random walk bridge $X_t, \, 0 \le t \le 2n$ such that 
\[  \Prob\{\sup_{0 \leq s \leq 1} |(2n)^{-1/2} X_{2ns} - B_{s}|
    \geq c \, n^{-1/2} \, \log n \}  \leq 
c \, n^{-30}.  \]
\end{lemma}
\begin{proof}  This is a special case of Theorem \ref{main}.
\end{proof}

\begin{remark}  The choice of $30$ as the exponent on the
right-hand side is arbitary.  The same
 result holds with error $O(n^{-r})$
for any $r > 0$, with suitably chosen $c = c_r$.
\end{remark}
\begin{corollary}  
There exists a $c < \infty$ such that  for every positive integer $n$,
 there exists a probability space $(\Omega,{\cal F},\Prob)$
on which are defined a two-dimensional
Brownian bridge $B_t, 0 \leq t \leq 1$
and a  two-dimensional simple
random walk bridge $X_t, \, 0 \le t \le 2n $ such that for
each $n$,
\[  \Prob\{\sup_{0 \leq s \leq 1} |n^{-1/2} X _{2ns} - B_{s}|
    \geq c \, n^{-1/2} \, \log n \}  \leq 
c \, n^{-30}. \]
\end{corollary}

\begin{proof}
 If   $S_j^{1},S_j^{2}$
are independent one-dimensional simple random walks, then
\begin{equation}  \label{twofromone}
       S_j := \frac{S_j^{1} + i \,  S_j^{2}}{1 + i} , 
\end{equation}
is a  two-dimensional simple random walk (written in
complex form).   Conditioning on   
  $S_{2n} = 0$  is the same as conditioning
on  $S_{2n}^1 = S_{2n}^2 = 0$.  In other words, we can
obtain a two-dimensional random walk bridge as the
product of two independent one-dimensional random
walk bridges.  Hence we can construct the probability space as
the product of two probability spaces as in the lemma.
\end{proof}
 The following 
corollary is in the form that we will need in the rest of the paper. We 
have chosen to write it out in detail because of this. Recall that if a 
process is defined only for integer times, we extend its definition to 
non-integer times by linear interpolation.

\begin{corollary} \label{jun2.cor1}
 There exists
a $c$ and a probability
space $(\Omega,{\cal F},\Prob)$ on which are defined
process $B^{n,z,m}, $ $ S^{n,z,m}$, $n=1,2,\ldots, z \in \Z^2,
m=1,2,\ldots$
such that
\begin{itemize}
\item  the processes 
\[ B^{n,z,m}_t, \;\;\;\;0 \leq t \leq 1;\;\;
 n = 0,1,\ldots; \;\; z \in \Z^2,\;\; m=1,2,\ldots; \]
 are independent two-dimensional Brownian bridges;
\item  the processes 
\[  S^{n,z,m}_j, \;\;\;    n = 0,1,\ldots; \;\; z \in \Z^2,
     j=0,\ldots,2n, \;\;m=1,2,\ldots ; \]
are independent and $S^{n,z,m}_j, j=0,\ldots,2n$ has the
distribution of a two-dimensional simple random walk conditioned
so that $S_{2n}^{n,z,m} = 0$.
\item \[ \Prob\{\sup_{0 \leq s \leq 1} |n^{-1/2} S _{2ns}^{n,z,m} - B_{s}^{n,z,m}|
    \geq c \, n^{-1/2} \, \log n  \}  \leq 
  c \, n^{-30}.  \]
\end{itemize}
\end{corollary}

\begin{proof}  Take products of the probability spaces
in the previous corollary.
\end{proof}

For the remainder of this paper we fix the probability space
$(\Omega ,{\cal F},\Prob)$ as in the previous corollary.  On this
probability space are defined the independent $\loops ^0$-valued
random variables
$\tilde L(n,z;m)$ as in Section~\ref{walksoup},  
\[ \tilde L(n,z;m) = [S^{n,z,m}_0,S^{n,z,m}_1,\ldots, S^{n,z,m}_{2n}] . \]
Also, we have independent, identically
distributed random variables $L(n,z;m)$
taking values in ${\cal C}_1$. 
\[ L(n,z;m)= B_{\cdot}^{n,z,m} \]
These have the distribution
of a two-dimensional Brownian bridge of time duration 1.
    We assume that these
are coupled as in the proposition. 

 \section{Constructing the Brownian  loop soup}

In this section we will show how to construct the Brownian loop soup
in a way that is highly correlated with the random walk loop soup.
We will restrict to the rooted Brownian loop soup restricted
to loops of time duration at least $5/8$; to get a complete realization
one can attach an independent realization of the loop soup
with loops of time duration less than $5/8$.  These small
loops will not be coupled with the random walk loops.
 Let $N(n,z;t)$ be a collection
of independent Poisson processes with parameter
\[ q_n :=          \int_{n-(3/8)}^{n+(5/8)} \frac{ds}{2\pi s^2} = \frac{1}{ 2\pi (n + (5/8))
 \, (n - (3/8))
     )}
  =  \frac{1}{2\pi n^2} - \frac 1 {8 \pi n^3}
                         +O(\frac{1}{n^4}). \]  

Recall the definition of $\tilde N(n,z;t)$ from Section~\ref{walksoup}, and note that 
$q_n - \tilde q_n = O(\frac{1}{n^4})$ .
(We have chosen to couple random walk loops of $2n$ steps with Brownian loops of
time duration $n - (3/8)$ to $n+ (5/8)$.  The particular choice of interval
$[-3/8,5/8]$ was chosen so that   $q_n$ and $\tilde q_n$ agree up to an error
of size $O(n^{-4})$.)
It is easy  to see that we can
  couple   $N(n,z;t), \tilde{N}(n,z,t)$ on the same probability space so that:
\begin{itemize}
\item  $\{N(n,z;t)\}$ are independent Poisson processes with parameter $q_n$;
\item  $\{\tilde N(n,z;t)\}$ are independent Poisson processes with
   parameter $\tilde q_n$ ;
\item There is a $c$ such that for all $n,z,t$,
$\Prob\{ N(n,z;t) \neq \tilde N(n,z;t)\} \leq t \,
|q_n -\tilde q_n| \leq ctn^{-4}.$
\end{itemize}
In fact, we can let $\hat N(n,z;t), n=1,2,\ldots, z\in \Z^2$, be independent
Poisson processes with parameter $1$ and then set
\[    N(n,z;t) = \hat N(n,z;q_nt), \;\;\;\;\;\;
              \tilde N(n,z;t) = \hat N(n,z,;\tilde q_nt) . \]
Assume without loss of generality that on this probability space we have independent copies of the coupled processes
\[                   L(n,z;m) , \;\; \tilde L(n,z;m)  \]
as in Section~\ref{strongapprox}, independent complex-valued
random variables
$Y(n,z;m)$ that are uniformly distributed on the square
$\{x+iy: |x| \leq 1/2, |y| \leq 1/2\}$, and independent real-valued random variables
$T(n,z;m)$  with density  
\[     \frac{(n+ \frac 58) (n - \frac 38)
 }{s^2},\;\;\;\;\;  n - \frac 38  \leq s \leq n+ 
  \frac 58 . \]                     

We construct the rooted Brownian loop soup (restricted to loops
of time duration at least $5/8$) as follows:
\begin{itemize}
\item $N(n,z;t)$ will be the number of rooted loops that have
appeared by time $t$ whose root is in the unit square centered at $z$ and
whose time duration is between $n-(3/8)$ and $n+(5/8)$; 
\item  scale the bridge (of time duration $1$ and rooted at $0$) $L(n,z,m)$
 so that it has time duration
$T(n,z,m)$; and then translate it so that its root is $z + Y(n,z,m)$;
we call this final loop $L^*(n,z,m)$.
\end{itemize}
Then it is easy to see from the definition
that the collection of loops
\[ {\cal A}_\lambda = \{L^*(n,z;m): N(n,z;\lambda) \geq m, \quad 
n \in \Z ^+, z \in \Z ^2 \} \]
is a realization of the Brownian loop soup with intensity $\lambda$ (restricted
to loops of time duration at least $5/8$).  We can then extend ${\cal A}_\lambda$
to a realization of the Brownian loop soup by adding an independent
realization of loops of time duration less than $5/8$.
Recall from the discussion in the Introduction that ${\cal A} _{\lambda,N}$ is 
also a realization of 
the Brownian loop soup.
On this space we also have the scaled random walk soup  $\tilde {\cal A} _{\lambda,N}$. 

We will now show that this coupling 
satisfies the conclusions of Theorem \ref{maintheorem2}.  Without
loss of generality, we may assume that $\lambda r^2 N^{2-3\theta}
  \leq 1$; in particular, $\lambda \leq N^4$.  First, note that
\begin{multline*} \hspace{.5in}\Prob\{N(n,z;\lambda) \neq \tilde N(n,z;\lambda)
 \mbox{ for some } N^{\theta} \leq n < \infty, |z| \leq r N \} \\
 \leq  c \sum_{|z| \leq r N} \sum_{n \geq N^{\theta}} \lambda \,
 n^{-4} \leq  c \, \lambda \, r^2 \, N^{2-3\theta}. \hspace{1.in}  
\end{multline*}
Hence, except for an event of probability $O(\lambda r^2 N^{2-3\theta})$,
$N(n,z;\lambda) = \tilde N(n,z;\lambda)$, for $n \ge N ^\theta, |z| \le rN$.
Also,
\[  \Prob\{N(n,z;\lambda) \geq N^5 \mbox { for some } N^{2/3} \leq n \leq N^6,
           |z| \leq r\, N\}  \leq  c \, r^2 \, N^8 \, \Prob\{Y \geq N^5\}
     \leq c \, r^2 \, N^{-5} , \]
where $Y$ is a Poisson random variable with expectation $c\, N^4$.  The last
estimate uses an easy estimate on Poisson random variables; in fact
for positive integers $N$,
$\Prob\{Y \geq N^5\} \leq \Prob\{Y \geq N^4\}^N \leq e^{-aN}.  $

Let 
\[     Z = Z_{N,r,\theta,\beta} =\sum_{|z| \leq rN} \sum_{n \geq N^{6}} 
                N(n,z,\lambda). \]
Then $Z$ is Poisson with 
\[ \E[Z] \leq c \,\lambda \, \sum_{|z| \leq rN} \sum_{n \geq N^{6}} 
  \frac{1}{n^2} \leq c \, \lambda \, r^2 \, N^{-4} . \]
Hence, $\Prob\{ Z \ne 0 \} \leq c \, \lambda\, r^2 N^{-4}\,  \leq
  c\, \lambda\, r^2 \, N^{2 - 3\theta}.$   Therefore,
\[ \Prob\{N(n,z,\lambda) >0 \mbox{ for some } n \geq N^{6},
  |z| \leq r N\} \leq c \, \lambda \, r^2 \, N^{2 - 3 \theta}  \]
and similarly the same estimate holds with $\tilde N(n,z,\lambda)$
replacing $N(n,z,\lambda)$.

Let us  denote the loops $L^*(n,z;m)$ and 
$\tilde L(n,z;m) + z $ by 
$\gamma_{n,z,m}$ and $\tilde\gamma_{n,z,m}$. 
Let $A = A_{N,r}$ be the event
\[  A = \{\sup_{0\le s \le 1} |\gamma_{n,z,m}(s t_\gamma)-\tilde\gamma_{n,z,m}(s t_{\tilde\gamma})|
 \ge c_2 \log N^6\hspace{2in}\]\[ \hspace{2in}
   \mbox{ for some } |z| < r N, N^{2/3} \leq n \leq N^6, m \leq N^5 \} . \]
Here we use $c_2$ for the constant $c$ from Corollary \ref{jun2.cor1}.  Then
the corollary tells us that
\[  \Prob(A) \leq r^2 \, N^2 \, N^6 \, N^5 \, O((N^{2/3})
      ^{-30}) \leq c\, r^2 \, N^{2 -
  3 \theta}.\]

On the intersection of $A^c$ and the events that
\[   N(n,z;\lambda) = \tilde N(n,z;\lambda) , \;\;\;\;  n \geq N^{\theta},
   \; |z| \leq r N , \]
\[  N(n,z;\lambda) \leq N^5, \;\;\;\; N^{2/3} \leq n \leq N^6,
 \; |z| < rN, \]
\[   N(n,z;\lambda) = \tilde N(n,z,\lambda) = 0, \;\;\;\; n \geq N^6,\;
   |z| \leq r N , \]
the coupling satisfies the conclusion of Theorem~\ref{maintheorem2}.

\section{Bounded domains}

In this section, we let $D$ denote a simply connected domain in $\C$
containing the origin, contained in the unit disk $\Disk = \{z \in \C:
|z| < 1\}$. If $\epsilon > 0$, let $D_\epsilon = \{z \in D:
\dist(z,\p D) > \epsilon\}$. Recall the definition of the loop measure
$\mu$ from the Introduction.  We say that a loop $\gamma$ is in
  $D$ if $\gamma[0,t_\gamma] \subset D$.

\begin{proposition} \label{disk}
There is a $c < \infty$ such that
if $0<\epsilon\leq t_0^{3/2}$, then
the $\mu$ measure of the set of
loops  in $\Disk$ of time duration at least
$t_0$ that are not in $\Disk_\epsilon$
 is bounded above by $c \, \epsilon \, t_0^{-3/2} $.
\end{proposition}

\begin{proof}   
The measure we are interested in is given by 
\[       \int_\Disk \int_{t_0}^\infty \frac{1}{2\pi t^2} 
\, p_t(z) \, dt \, dz, \]
where $p_t(z) =  \Prob\{0 < \dist(B^{\rm br}_z\,[0,t],\p \Disk) \leq \epsilon\}$ 
and $B^{\rm br}_z$ denotes a Brownian bridge of time duration $t$ rooted 
at $z$.  By time reversal,
we can see that, 
 \[ \Prob\{0 < \dist(B^{\rm br}_z\,[0,t],\p \Disk) \le \epsilon\} \le
2\, \Prob\{0 < \dist(B^{\rm br}_z\,[0,\frac{t}{2}] ,\p \Disk) \le \epsilon, 
B^{\rm br}_z\,[0,t] \subset \Disk \}.
\]
But,  
\begin{multline*}
\Prob\{0 < \dist(B^{\rm br}_z\,[0,\frac{t}{2}] ,\p \Disk) \le \epsilon, 
B^{\rm br}_z\,[0,t] \subset \Disk \} = \\
\lim _{\delta \to 0^+} \frac{2t}{\delta^2}
\Prob^z\{0 < \dist(B[0,\frac{t}{2}] ,\p \Disk) \le \epsilon, 
B[0,t] \subset \Disk ,|B_t -z| < \delta\},
\end{multline*}
where $B_t$ denotes   a standard Brownian motion.
By the strong Markov property and the ``gambler's ruin'' 
estimate\footnote{The gambler's ruin estimate  
 states that the probability that a one dimensional
standard Brownian
motion starting at $\epsilon > 0$ stays positive up to time $t$
is bounded above by $c \, \epsilon \, t^{-1/2}$.}, 
\[  \Prob^z\{0 < \dist(B[0,\frac{t}{2}],\p \Disk) \leq \epsilon;
               B[0,\frac 34 t] \subset \Disk\} \leq c \, \epsilon \,
            t^{-1/2} . \]
Given this event, the probability that $|B_t - z| < \delta$ is bounded
above by $c \delta^2/t$.  Hence
\[          p_t(z) \le c \,\epsilon \,\delta\,  t^{-1/2} , \]
and the result follows by integrating.
\end{proof}

\begin{proposition} \label{general}Suppose $D$ is a simply
connected domain contained in the unit disk.
There is a $c < \infty$ such that
if $0< \epsilon \leq t_0^{5/4}$, then
the $\mu$ measure of the set of
loops  in $D$ of time duration at least
$t_0$ that are not in $D_\epsilon$
 is bounded above by $c \, \epsilon^{1/2} \, t_0^{-5/4} $.
\end{proposition}

\begin{proof}  
The measure we are interested in is  
\[       \int_D \int_{t_0}^\infty \frac{1}{2\pi t^2}\,  p_t(z) \, dt \, dz, \]
where $p_t(z) =  \Prob\{0 < \dist(B^{\rm br}_z\,[0,t],\p D) \leq \epsilon\}$ 
and $B^{\rm br}_z$ denotes a Brownian bridge of time duration $t$ rooted 
at $z$.  By time reversal,
we can see that, 
 \[ \Prob\{0 < \dist(B^{\rm br}_z\,[0,t],\p D) \leq \epsilon\} \le
2 \Prob\{0 < \dist(B^{\rm br}_z\,[0,\frac{t}{2}] ,\p D) \le \epsilon, 
B^{\rm br}_z\,[0,t] \subset D \}.
\]
But, 
\begin{multline*}
\Prob\{0 < \dist(B^{\rm br}_z\,[0,\frac{t}{2}] ,\p D) \leq  \epsilon, 
B^{\rm br}_z\,[0,t] \subset D \} = \\
\lim _{\delta \to 0^+} \frac{2t}{\delta^2}\,
\Prob^z\{0 < \dist(B[0,\frac{t}{2}] ,\p D) \leq  \epsilon, 
B[0,t] \subset D ,|B_t -z| < \delta\},
\end{multline*}
where $B_t$ denotes   a standard Brownian motion.
By the strong Markov property and the Beurling estimate 
(see Lemma~\ref{beurling} below),
\[  \Prob^z\{0 < \dist(B[0,\frac{t}{2}],\p D) \leq \epsilon;
               B[0,\frac 34 t] \subset D\} \leq c \, \epsilon^{1/2} \,
            t^{-1/4} . \]
Given this event, the probability that $|B_t - z| < \delta$ is bounded
above by $c \delta^2 /t$.  Hence
\[          p_t(z) \leq c \,\epsilon^{1/2} \,  t^{-1/4} , \]
and the result follows by integrating.
\end{proof}

\begin{lemma}[Beurlng estimate] \label{beurling} Let $B_t$ denote a standard
two-dimensional Brownian
motion. There is a $c$ such that if  
$\gamma:[0,\infty) \to \C$ is any 
continuous curve with $|\gamma(0)| = r$ 
and $\lim _{t \to \infty} |\gamma(t)| = \infty$, then 
\[ \Prob ^z \{B[0,t] \cap \gamma[0,\infty) = \emptyset \} 
\le c \left(\frac{r}{\sqrt{t}}\right) ^{\frac{1}{2}}, \qquad   
  |z| \le r.
\]
\end{lemma}
\begin{proof}   
If $t$ is replaced by $\sigma_t = \inf\{s: |B_s| = \sqrt{t}\}$,
then this lemma follows from the Beurling estimate which is a
corollary of the Beurling Projection Theorem
\cite[Theorem V.4.1]{beurling.est}.  We will assume that estimate and show how
the estimate  for fixed times $t$ can be
deduced from the result
for the stopping times $\sigma_t$. 
By scaling it is enough to do $t=1$.  Let 
$\tau _n = \inf \{t : |B_t - B_0| = 2 ^ {-n} \}.$
We will show that $\Prob \{ B[0,1\wedge \tau _0] \cap \gamma[0,\infty) 
= \emptyset \} \le c r ^{\frac{1}{2}}.$

First, note that 
\begin{multline*}
  \Prob ^z\{ B[0,1\wedge \tau _0] \cap \gamma[0,\infty) = \emptyset \}  \\
 \le \Prob ^z\{ B[0, \tau _0] \cap \gamma[0,\infty) = \emptyset \} 
 + \Prob ^z\{B[0, 1] 
\cap \gamma[0,\infty) = \emptyset, \tau _0 > 1 \}.
\end{multline*}
The Beurling estimate  gives us that the first 
summand is bounded by $cr^{1/2}$. Thus, we only need to bound 
the second summand.
To do this we let $A_k$ be the event
 \[ A_k = \{ \tau _k - \tau _{k+1} > 3 ^{-(k+1)}
 \text{ and } 
\tau _j - \tau _{j+ 1} \le 3 ^{-(j+1)}, \text{ for all $j > k$ }
\}. \]

Note that
\begin{equation}  \label{jose}
  \Prob^z\{\tau _k - \tau _{k+1} > 3 ^{-(k+1)}\}
   = \Prob^0\{\tau_0 - \tau_1 > (4/3)^k/3 \} \leq c_1
  \exp\{-c_2 \, (4/3)^k \} . 
\end{equation}
The last estimate is a standard estimate for Brownian motion
and implies that $\Prob\{\tau_k - \tau_{k+1} > 3^{-(k+1)}
\mbox{ i.o. } \} = 0 .$
Since 
 $\sum _{k=1} ^\infty 3 ^{-k} < 1$, we see that
the event $\{\tau_0 > 1\}$
is contained in $\cup_{k \geq 0} A_k$ up to an event
of probability zero.
 Using  the Beurling estimate, 
 we see that  
\begin{eqnarray*}
\lefteqn{\Prob ^z \{ B[0, 1] \cap \gamma[0,\infty) = \emptyset, \tau _0 > 1 \}}
  \hspace{1.3in}\\ &\leq&
\sum _{k = 0 } ^ \infty \Prob ^z \{ B[0, 1] \cap \gamma[0,\infty) = \emptyset,
A_k \} \\
& \leq & \sum _{k = 0 } ^ \infty \Prob^z\{  B[0,\tau_{k+1}] \cap \gamma[0,
\infty) = \emptyset, \tau_k - \tau_{k+1} \geq 3^{-(k+1)}\}\\
&\le & \sum _{k = 0 } ^ \infty c \left(\frac {r} {2 ^{-(k+1)}} \right) ^{1/2}
\Prob^z \{ \tau _k - \tau _{k+1} > 3 ^{-(k+1)} \}  \le cr^{1/2}.
\end{eqnarray*}
 The last inequality uses (\ref{jose}).

\end{proof}

\begin{corollary} There is a $c$ such that for every $N, \lambda$ and 
$\theta < 2$; there exists a coupling of the Brownian loop soup 
 restricted to
$\Disk\,[D]$ and a (1/N)-random walk soup restricted to 
$\Disk\,[D]$
such that the one-to-one correspondence of Theorem~\ref{maintheorem2} holds 
except on an event of probability at most 
$c(\lambda + 1) \log N \,N^{2 - (2/3)\theta}$ 
[respectively, $\,c(\lambda + 1) \log N^{(1/2)}\, N^{2-(5/4)\theta}]$.
\end{corollary}
\begin{proof}
It follows from Proposition~\ref{disk} that the probability that a realization
of the Brownian loop soup has at least one 
loop of time duration greater than $cN^{\theta -2}$ 
staying in $\Disk$, but not 
in $\Disk_\epsilon$ for $\epsilon = c(\log N/\,N)$ is 
$O(\lambda \log N \, N^{2-(3/2)\theta})$. For general domains we get
$O(\lambda (\log N)^{(1/2)} \, N^{2-(5/4)\theta})$ upon using Proposition~\ref{general}. Therefore, if we consider the coupling of a $(1/N)$-random
walk soup and a Brownian soup and we restrict to those loops 
in a domain $D$, then we get the one-to-one correspondence of the loops
as before, except on an event of probability 
$c \,(\lambda + 1) \, N^{2 - (3/2) \theta}
 \, \log N$ if $D = \Disk$ or an event of probability
$c \,( \lambda + 1) \, N^{2 - (5/4)\theta} \, \log N^{1/2}$ for general 
simply connected $D$ contained in the unit disk.  
  
\end{proof}

\section{The dyadic approximation}  \label{kmtsec}

\subsection{Introduction}

In this note we give a proof of the ``dyadic''  
   strong approximation for random walk bridges
 by Brownian bridges using the methods in \cite{KMT}.
Let $X_1,X_2,\ldots$ be independent random variables
with $\Prob\{X_j = 1 \} = \Prob\{X_j = - 1\} = 1/2$
and let $S_n = X_1 + \cdots + X_n$. For  positive integers
$n$, let
$L_n = \{z \in \Z: \Prob \{S _n = z \} > 0 \}$.
If $z \in L_n$, $\{S^{(n,z)} _m\}_{m=0}^n$ will denote
 a process with the law 
of $\{S_m\}_{m=0}^n$ conditioned so that $S_n = z$. 

  We start with
a definition.  Suppose $Z$ is a continuous
random variable with strictly
increasing distribution function $F$ and $G$ is the distribution
function of a discrete random variable whose support
is $\{a_1,a_2,\ldots\}$.  Then $(Z,W)$ are
{\em quantile-coupled (with distribution functions $(F,G)$)}
 if $W$ is defined by
\[    W = a_j \;\;\;\mbox{ if } \;\;\; r_{j-} < Z \leq r_j , \]
where $r_{j-},r_j$ are defined by
\[     F(r_{j-}) = G(a_j-), \;\;\;\;\;
           F(r_{j}) = G(a_j). \]
The quantile-coupling has the following property. If  
\[  F(a_k - x) \leq G(a_{k} -) < G(a_k) \leq F(a_k + x) , \]
then
\begin{equation}  \label{may29.1}
   |Z-W| = |Z-a_k|
 \leq x \;\;\;\;\mbox{on the event} \;\;\;\;
     \{W = a_k\}. 
\end{equation}
We will need the following lemmas about the random walk;
we will prove them in \S\ref{lemmasec}.

\begin{lemma}  \label{stephen.1}
There exists $\epsilon_0 > 0$ such that for
every $b_1 > 0$ there exist $0 < c,a < \infty $ such that
the following holds.  Let $N$ be a $N(0,1)$ random variable.
For each integer $n > 1$, each  integer
 $m  $ with $|2m-n| \leq 1$,
and every $z \in L_n$, let
\[    Z = Z^{(m,n,z)} = \frac mn z + \sqrt{m (1 - \frac mn)} 
\; N, \]
so that $Z \sim N(\frac mn z, m(1-\frac mn))$.
Let $W = W^{(m,n,z)}$ be the random variable with
${ \cal L } (W) = { \cal L} (S^{(n,z)} _m)$ that is
quantile-coupled with $Z$.  Then if $|z| \leq \epsilon_0n$ and
$\Prob\{W = w\} > 0$, 
\begin{equation}
\label{stephen.5}
   \E[e^{a| Z-W|} \mid W= w ]\leq c \, \exp\left\{b_1
   \frac{w^2 + z^2}{n} \right\} . 
\end{equation}
\end{lemma}

\begin{remark}  For simple random walk, it is easy to show
that (\ref{stephen.5}) holds for $ \epsilon_0 n
\leq |z| \leq n$ (with perhaps
different $a,c$), so it follows that the result holds for all
$|z|$.  However, we state the lemma only for $|z| \leq \epsilon_0n$
 because this is all
that we use.   
\end{remark} 

\begin{lemma}  \label{stephen.2}
There exist $c_2, b_2,\epsilon_0$ such that 
for every integer
$n \ge 2$, every integer $m$ with $|2m-n| \le 1$, 
every $z \in L_n$ with $|z| \leq \epsilon_0 n$,   and every
$w \in \Z$,
\[  \Prob\{S_m = w \mid S_n = z\}
              \leq c_2 \, n^{-1/2} \, \exp\left\{-b_2 \frac{(w -
   (z/2))^2}{n}\right\}  . \]
\end{lemma}
 
\begin{remark}  We can actually
show that this holds for any
$\epsilon_0 < 1$ (with the constants $c_2,b_2$ depending
on $\epsilon_0$).
\end{remark}

   Let $B$ denote
a  Brownian bridge, i.e., a Brownian motion
in $\R$ conditioned so that $B_0 = B_1 = 0$ (see \S\ref{bridgesec}
 for a more precise definition).  If $z_1,z_2 \in \R$,
$n > 0$,  then
\begin{equation}  \label{sept15}
   Y_t^{(n,z_1,z_2)} := \sqrt n \, B_{t/n} + \frac{n-t}{n}\, z_1 + \frac{t}{n} \, z_2 , \;\;\;\;
   0 \leq t \leq n , 
\end{equation}
is the Brownian bridge conditioned so that $B_0 = z_1,
B_n = z_2$.  We write $Y_t^{(n,z)}$ for
$Y_t^{(n,0,z)}$.   If $(S ^{(n,z)},B)$ are defined on the same probability
space, we define
\[  \Delta(n,z) = \Delta(n,z,S^{(n,z)},B)
    = \sup_{0 \leq t \leq n} |Y_t^{(n,z)} - S^{(n,z)}_t| . \]
(Recall that  $S^{(n,z)}_t$ is defined for noninteger $t$ by
linear interpolation.)  In Section \ref{mainproofsec}, we
will prove the following.

\begin{theorem} \label{maintheorem} For every $b > 0$,
there exist $0 < c,a,\alpha < \infty$
such that for every positive integer $n$, 
there is a probability space on which are defined 
a Brownian bridge $B$ and the family of processes 
$\{S ^{(n,z)}: z \in L _n \}$ such that
\begin{equation}  \label{theoremeq}
\E[e^{a \Delta(n,z)}] \leq
               c \, n^{\alpha} \, e^{b|z|^2/n} . 
\end{equation}
\end{theorem}

Using Chebyshev's inequality, we get the following theorem
as a corollary.

\begin{theorem}  \label{main}
For
every $b > 0$  there exist $0 < c,\alpha < \infty$, such that
for every postive integer $n$, 
there is a probability space on which are defined 
a Brownian bridge $B$ and the family of processes 
$\{S ^{(n,z)}: z \in L _n \}$ such that for all  
$r > 0$,
\[   \Prob\{\Delta(n,z) > r \, c \, \log n \}
   \leq c \, n^{\alpha - r}\, e^{bz^2/n} . \]
\end{theorem}

\subsection{Brownian bridge}  \label{bridgesec}

If $W_t$ denotes a standard one-dimensional Brownian motion, then the
process
\[  B_t = W_t - t \, W_1, \;\;\;\; 0\leq t \leq 1 , \]
is called a {\em Brownian bridge (conditioned so
that $B_0 = 0, B_1 = 0$)}.  It can also be characterized
as the continuous Gaussian process $B_t, 0 \leq t \leq 1$
with
\[      \E[B_t] = 0 \qquad
           {\bf Cov}[B_s,B_t]=\E[B_sB_t] = s (1-t) , \;\;\; 0 \leq s \leq t
   \leq 1 . \]
 More generally, if
$B_t$ is a Brownian bridge and
\[          X_t = \sqrt{s_2 - s_1} \,
             B_{(t-s_1)/(s_2-s_1)} \; +
         x_1 + (\frac{t - s_1}{s_2 - s_1}) \, (x_2 - x_1), \]
is the Brownian bridge conditioned so that
$X_{s_1} = x_1, X_{s_2} = x_2$.  It is the continuous
Gaussian process $X_t, s_1 \leq t \leq s_2$ with
\[ \E[X_t] = x_1 + (\frac{t - s_1}{s_2 - s_1}) \, (x_2 - x_1)
      , \;\;\;\;
  {\bf Cov}[X_sX_t] =     \frac{(s - s_1) \; (s_2 - t)}{s_2 - s_1}  \, \;\;\;\;
   s_1 \leq s \leq t \leq s_2 . \]

\begin{lemma} \label{surgery}
Suppose $B,\tilde B$ are independent Brownian
bridges and $N$ is an indepedent $N(0,1)$ random variable.
Suppose $0 < s < 1$, and  define
$X_t, 0 \leq t \leq 1$ by
\[   X_s =
 \sqrt{s(1-s)} \, N \]
\[      X_t = \sqrt{s} \,
             B_{t/s} \;  
         + \frac ts \, X_s, \;\;\;
            0 \leq t \leq s , \]
\[    X_t =  \sqrt{1 - s} \,
             \tilde B_{(t - s)/(1 - s)} \; +
        \frac{1-t}{1-s}\,  X_s  , \;\;\;
      s \leq t \leq 1. \]
Then $X_t$ is a Brownian bridge conditioned so that
$X_0 = X_1=0  . $
\end{lemma}

\begin{proof} This can be easily checked using the
Gaussian characterization of Brownian bridges.  The formulas
are not mysterious.  What we are doing is defining 
$X_t, 0 \leq t \leq 1,$ by first choosing $X_s$
(using the appropriate distribution on $X_s$), then defining
$X_t, 0 \leq t \leq s$, and $X_t, s \leq t \leq 1$,
as appropriate Brownian bridges. 
\end{proof}

\medskip  
We will need the following easy estimate for Brownian bridges.
Let
\[  M  = \sup_{0 \leq t \leq 1} |B_t| . \]
Then there exist $\tilde c,u$ such
that for all $a > 0$,
\begin{equation}  \label{1.1}
              \E[e^{a M}] \leq \tilde c \, e^{ua^2} . 
\end{equation}
If $B_t$ is replaced by a Brownian motion $W_t$, this
estimate is standard using the reflection principle.
That argument can easily be adapted to establish (\ref{1.1}),
perhaps with different $\tilde c, u$. (In fact, the maximum
for Brownian motion stochastically dominates $M$ so (\ref{1.1})
  holds with the same $\tilde c, u$, but we will not
need this stronger fact.)

\subsection{Proof of Theorem \ref{maintheorem}}  \label{mainproofsec}

It suffices to prove the result for $b$ sufficiently small.  
We fix positive $b < b_2/37$ where $b_2$ is the constant
from Lemma \ref{stephen.2}.  We let $\epsilon_0$ be the
smaller of the two values of $\epsilon_0$ in Lemmas 
\ref{stephen.1} and \ref{stephen.2}.

In this proof, by an {\em $n$-coupling} we will mean
a probability space on which are defined 
a Brownian bridge $B$ and the family of processes 
$\{S ^{(n,z)}: z \in L _n \}$.

Note that for any $n$-coupling, if $z \in L_n$,
 $S_t = S_t^{(n,z)},$ and
$ Y_t = Y_t^{(n,z)}$ as in (\ref{sept15}),
then
\[   \Delta(n,z) =
         \sup_{0 \leq t \leq n} |S_t - Y_t|
     \leq \sup_{0 \leq t \leq n} |S_t| + \sup_{0 \leq t
  \leq n} |Y_t|
      \leq   n + \left[n + \sqrt n \sup_{0 \leq t \leq 1} |B_t| 
\right]. \]
Hence,
\[   \E[e^{a \Delta(n,z)}
  ]   \leq e^{2an} \, \E[\exp\{ a \sqrt n 
                   [\sup_{0 \leq t \leq 1} |B_t|] \}]
    \leq \tilde c\, e^{(2a + u a^2) n} , \]
where $u,\tilde c$ are as in  (\ref{1.1}). Clearly, there exists 
$a_0 =a_0(b) >0$ such that if $a \in (0,a_0)$, then 
$2a + u a^2 \leq b\, \epsilon_0^2$.

Therefore, for any $n$-coupling 
inequality (\ref{theoremeq}) will hold with 
$c = \tilde c, \alpha = 0$ and $a \in (0,a_0)$
for all $z \in L_n$ with $|z| \ge n\epsilon_0$.
For the remainder of this section, we will assume
$a < a_0$.
We wil also assume that $a$ is sufficiently small so that
(\ref{stephen.5}) holds with $b_1 = b/20$.  We now fix
such a value of $a$, and we will show how to construct the
$n$-couplings so that (\ref{theoremeq}) holds
for some $c,\alpha$.

We will use an induction.
Clearly, we can choose $n$-couplings for $n \le 2$ such that
\[
\E[e^{a \Delta(n,z)}] e^{-b|z|^2/n} \leq C \qquad \forall z \in L_n,\, n \le 2.
\]
We can assume without loss of generality that $C \ge 1$.
We will show that there exists a constant $c$ (which
without loss of 
generality   we can assume is
greater than both 1 and the $\tilde c$ 
in (\ref{1.1})) such that for  every
positive integer 
$s$, if there exist $n$-couplings for all $n \le 2^s$ such that 
\begin{equation}\label{preclaim}
\E[e^{a \Delta(n,z)}] e^{-b|z|^2/n} \leq C(s),
\end{equation}
then there exist $n$-couplings for all $n \le 2^{s+1}$ such that 
\begin{equation}\label{claim}
\E[e^{a \Delta(n,z)}] e^{-b|z|^2/n} \leq c\,C(s). 
\end{equation}
The theorem follows easily from this claim.
 
In order to prove the claim above, let $2^s < n \leq 2^{s+1}$.
We will show how to construct  a probability space
on which are defined  a Brownian 
bridge and a family of processes $\{S ^{(n,z)}: z \in L _n, \, |z| \le n\epsilon_0 \}$ 
satisfying~(\ref{claim}). Once this is done, we can adjoin, possibly after
enlarging the probability space, the processes for $|z| > n\epsilon_0$. 
Since
  $c \ge \tilde{c}$, (\ref{claim}) will hold for these processes also. 
Hence, we assume
$|z| \leq \epsilon_0 n$. For notational ease we
will assume that $n$ is even
and we write  $n = 2k$. Note that
$k$ is an integer  with $2^{s-1} < k \leq 2^s$. 
(If $n$ is odd we write $n = k + (k+1)$ and do a similar
argument.)

We define the $n$-coupling as follows:
\begin{itemize}
\item Choose two independent $k$-couplings 
\[\left(\{ S ^{1\,(k,z)} \} _{z \in L_k}, B^1 \right), \qquad
\left(\{ S ^{2\,(k,z)} \} _{z \in L_k}, B^2 \right),
\]
satisfying~(\ref{preclaim}).
\item Let $N \sim N(0, 1)$ and define the translated 
{\em normal} random variables $Z^z = \sqrt{n/4} \, N + \frac{z}{2}$. Define the 
quantile-coupled random variables $W^z$ as in Lemma~\ref{stephen.1}. Assume, 
as we may, that all these random variables are independent of the 
two $k$-couplings chosen above. Note that $a$ has been chosen sufficiently small
so that (\ref{stephen.5}) holds with $b_1 = b/20$;
i.e.,
\[   \E[e^{a |Z^z - W^z|} \mid W^z = w]
 \leq c \, \exp[\frac{b}{20} \, \frac{w^2 + z^2}{n}] .
\]
\item Let 
\begin{equation}
  B_t =
  \begin{cases}
    \frac{1}{\sqrt{2}}B^1 _{2t} + t\, N& 0 \le t \le \frac{1}{2}, \\
    \frac{1}{\sqrt{2}}B^2 _{2(t -\frac{1}{2})} + (1-t)\, N& \frac{1}{2} \le t \le 1.
  \end{cases}
\end{equation}
By Lemma \ref{surgery}, $B_t$ is a Brownian bridge.
\item Let $S_k^{(n,z)} = W^z$, and
\begin{equation*}
S ^{(n,z)} _m = 
\begin{cases}
  S ^{1\,(k,W^z)} _m& 0 \le m \le k, \\
W^z + S ^{2\,(k,z-W^z)} _{m-k}& k \le m \le n. 
\end{cases}
\end{equation*}
What we have done is to first choose the value of $S^{(n,z)}_k$
from the appropriate conditional distribution of $S_k$ given
$S_n = z$  and then to define
the other values of $S^{(n,z)}_m$ from the conditional distribution
of $S_m$ given $S_n = z, S_k = W^z$.  
\end{itemize}
This gives us our coupling; we need to show that it
satisfies (\ref{claim}).
Note that 
\[  \Delta(n,z,S^{(n,z)},B)  \leq \hspace{3.5in}\]\[|Z^z - W^z| + 
  \max\left\{\, \Delta(k,W^z,S^{1\,(k, W^z)},B^1) \, ,\, 
            \Delta(k,z-W^z,S^{2\,(k,z-W^z)},B^2)\, \right\} . \]
Therefore,
\[  \E[e^{a \Delta(n,z)} \mid
W^z = w] \leq 
  \E[e^{a |Z^z - W^z|} \mid W^z = w]
\; C(s)\;  (e^{b|w|^2/k} + e^{b|z-w|^2/k})  . \]
Here we have used the
the fact that our two $k$-couplings satisfy
(\ref{preclaim}) and the  simple inequality,
$\E[e^{\max\{Z_1,Z_2\}}] \leq \E[e^{Z_1}] +
 \E[e^{Z_2}]$.
    Therefore,
\begin{equation*}
 \frac{\E[e^{a \Delta(n,z)}]}{C(s)}  \leq  c \, \sum_{w} \Prob\{W^z=w\}\;
          \exp\left\{\frac b {20} \, \frac{w^2 + z^2}{n}\right\}
                \exp\left\{b\frac{\max\{w^2, (z-w)^2\}}{k}\right\},
\end{equation*}
Since $w^2 + z^2 \leq 5 \, \max\{w^2, (z-w)^2\}$ and $k = n/2$,
this sum is bounded by
\[ c\, \sum_{w} 
\Prob\{W^z=w\}\;
  \exp\left\{\frac 9 4 \, b\, \frac{\max\{w^2, (z-w)^2\}}
               {n} \right\} . \]
We now split this sum into two pieces:
$|w-\frac{z}{2}| \leq |z|/6$ and $|w- \frac{z}{2}| > |z|/6$.  If
$|w-\frac{z}{2}| \leq |z|/6$, then $\max\{w^2,(z-w)^2\} \leq
(2z/3)^2$; hence,
\[ c\, \sum_{|w - \frac{z}{2}| \leq |z|/6}
\Prob\{W^z=w\}\;
  \exp\left\{\frac 9 4 \, b\, \frac{\max\{w^2, (z-w)^2\}}
               {n} \right\}  \leq c \, \exp\{\frac{bz^2}{n}\} . \]
To handle the other piece we make use of Lemma \ref{stephen.2}. Since
 $b_2 >37 b$, we have 
\[  \Prob\{W^z=w\}=\Prob\{S_k = w \mid S_n = z\}
   \leq c \, n^{-1/2} \;  \exp\{-37 b \, \frac{(w-\frac{z}{2})^2}{ n}\}
                    ,\]
and so, the sum over the second piece is bounded by
\[ c \, \sum_{|w-\frac{z}{2}| > |z|/6}  
   n^{-1/2} \;  \exp\{-37 b \, \frac{(w-\frac{z}{2})^2}{ n}\}
  \; \exp\{\frac 9 4 \, b\, \frac{\max\{w^2, (z-w)^2\}}
               {n} \} 
. \]
 For $|w-\frac{z}{2}| > |z|/6$, we have $(w-\frac{z}{2})^2 > \frac 1 {16} \, \max\{w^2, (z-w)^2\}.$
Hence,   the sum is bounded by
\[ c \, \sum_{w}  
   n^{-1/2} \;  \exp\{- \frac b {16}  \, \frac{w^2}{ n}\}
, \]
which is clearly bounded by a constant. 
 Therefore, 
\[
\frac{\E[e^{a \Delta(n,z)}]}{C(s)} \le c \, \exp\{\frac{bz^2}{n}\}.
\]
That is, our $n$-coupling satisfies equation~(\ref{claim}).
 
\subsection{ Local central limit theorem}
  \label{lemmasec}

We will derive the local central limit theorem for conditioned
simple random walks which is essentially a normal approximation
for hypergeometric random variables.  Although this is standard,
we derive it here because the size of the error term is important
for us.  Our starting point is Stirling's formula with
error
\begin{equation}  \label{stirling1}
       n! =  \sqrt{2 \pi} \, n^{n+(1/2)} \, e^{-n} \left[
1 + O(n^{-1})\right].
\end{equation}

\begin{lemma} 
Suppose $l,m,j$ are integers with $m >0, |l| \leq m/2, |j| \leq m/8$.
  Then
\[  \Prob\{S_{2m} =2j + 2l \mid S_{4m} = 4l \}
   = \hspace{2.8in} \]
\[2 \, \sqrt {\frac{1}{2\pi\,m(1-(l/m)^2)}}\;
   \exp\left\{-\frac{ (2j)^2}{2\,m(1-(l/m)^2)} + O(\frac{1}{m}) +
     O(\frac{j^4 }{m^3}) \right\}. 
\]
\end{lemma}

\begin{proof}
  Throughout this proof we assume that $|l| \leq m/2,
|j| \leq m/8$.  Without loss of generality we may
assume that $l \geq 0$. Let
\[  p(2m,2j\mid 4l) = \Prob\{S_{2m} = 2j \mid
     S_{4m} = 4l \}  = \frac{\binom{2m}{m +j}
  \binom{2m}{m + 2l - j}}{\binom{4m}{2m+2l}}. \]
Then,
\[ p(2m,2l + 2(j+1) \mid 4l) =
             p(2m, 2l + 2j \mid 4l) \; \frac{(m-j)^2 - l^2}
   {(m+j+1)^2 - l^2} , \]
and if $j > 0$,
\begin{eqnarray*}
 p(2m,2l + 2j \mid 4l) & =  & p(2m,2l \mid 4l) \; \prod_{i=0}
^{j-1} \frac{(m-i)^2 - l^2}
   {(m+i+1)^2 - l^2} \\
 & =  & p(2m,2l \mid 4l) \;   \prod_{i=1}
^{j} [1 - \frac{4im + 2i - 2m -1}{(m+i)^2 - l^2}] .
\end{eqnarray*}
By symmetry, $p(2m,2l-2j \mid 4l) = p(2m,2l + 2j \mid 4l)$.

Stirling's formula  (\ref{stirling1})
gives
\begin{equation}  \label{jan2.4}
p(2m,2l \mid 4l)= \frac{{2m \choose m+l}^2}
     {{4m \choose 2m+2l}} =2 \, \sqrt {\frac{1}{2\pi\,m(1-(l/m)^2)}}
\; \left[1 + O(\frac 1m)\right]
\end{equation}
This gives the result for $j=0$. To finish the proof of the result we need to 
show that 
\[ \sum _{i=1} ^j \log 
[1 - \frac{4im + 2i - 2m -1}{(m+i)^2 - l^2}]
= -\frac{ (2j)^2}{2\,m(1-(l/m)^2)} + O(\frac{1}{m}) +
     O(\frac{j^4 }{m^3}).
\]

Note that if $l \leq m/2$, $ j \leq       m/8$ and
$1 \leq i \leq j$, then
\[     
 0 \le \frac{4im + 2i - 2m -1}{(m+i)^2 - l^2} \leq \frac23. \]
There exists $c_1$ such that
\[  |\log (1-x) + x + \frac{x^2}{2}| \leq c_1 x^3, \;\;\;\;
   0 \leq x \leq 2/3. \]
Therefore,
\[ \sum _{i=1} ^j \log 
[1 - \frac{4im + 2i - 2m -1}{(m+i)^2 - l^2}]
  = \hspace{2in} \]
\[  O(\frac{j^4}{m^3}) -
\sum_{i=1}^j  \frac{4im + 2i - 2m -1}{(m+i)^2 - l^2}
   -\frac 12 \sum_{i=1}^j  (\frac{4im + 2i - 2m -1}{(m+i)^2 - l^2}
)^{2}  . \]
Note that
\[  \frac{1}{(m+i)^2 - l^2} = \frac{1}{m^2 - l^2} 
   -\frac{2im}{(m^2 - l^2)^2} + O(\frac{i^2}
   {m^4}).  \]
Hence,
\begin{eqnarray*} 
{\sum_{i=1}^j  \frac{4im + 2i - 2m -1}{(m+i)^2 - l^2}}
            & = & O(\frac 1m +\frac{j^4}{m^3}) +
 \sum_{i=1}^{j}  \frac{4im - 2m}{m^2 - l^2}
  - \sum_{i=1}^j \frac{8i^2m^2}{(m^2 - l^2)^2} \\
& = & O(\frac 1m + \frac{j^4}{m^3})+ \frac{2j(j+1)m - 2jm}
                    {m^2 - l^2} - \frac{(8/3) j^3m^2}{(
  m^2 -l^2)^2}
 \\
& = & O(\frac 1m + \frac{j^4}{m^3}) + 
      \frac{2j^2m}{m^2 -l^2} - \frac{(8/3) j^3m^2}{(
  m^2 -l^2)^2}  
\end{eqnarray*}
Also,
\[ \frac 12 \sum_{i=1}^j  (\frac{4im + 2i - 2m -1}{(m+i)^2 - l^2})^2
  = O(\frac 1m + \frac{j^4}{m^3}) + \frac{(8/3) j^3m^2}{(
  m^2 -l^2)^2}. \]
The result follows immediately.
\end{proof}

One can deduce the following more general case from the result above.
\begin{lemma}\label{clclt}
There is a $c$ such that for $n \ge 2$ an integer,
 and $m$ an integer with $|2m-n| \le 1$. 
If $|z|,|w| \le c/n;\,\, z \in L_n$ and $w+\frac mn z \in L_m$. Then 
\begin{align*}
\Prob \{ S _m = \frac mn z + w | S _n = z \} =& \\
2 \sqrt{ \frac {1}{2 \pi (m - (m^2/n)) (1 - (z/n)^2)}} 
&\exp \left\{ -\frac{w^2}{2 (m - (m^2/n)) (1-(z/n)^2)} + 
O( \frac 1n  + \frac {w^4}{n^3})\right\}.
\end{align*}
\end{lemma}

 We now state without
proof an easy large deviation estimate that follows from
large deviations for binomial random variables.

\begin{lemma} \label{largedev}
There exists an $\eta > 0$ such that, for any $a>0$, 
there exist $C = C(a) <\infty$, and
 $\gamma=\gamma(a)> 0$, such that 
for all $z$ with $|z|/n < \eta$ 
\[
\Prob \{ |S_m - \frac{m}{n} z | > am | S _n =z \} \le Ce ^{-\gamma m}
\]
\end{lemma}

 Lemma \ref{stephen.2} follows easily from Lemmas \ref{clclt} and 
\ref{largedev}.

\subsection{Coupling of conditioned distribution and normal}
\label{lemmasec2}

In the remainder of this section we prove Lemma~\ref{stephen.1}. Note that 
we only need to prove the lemma for $n$ sufficiently large. In order 
to simplify the notation we will assume that $n$ is even and hence $m=n/2$.
If $n$ is odd, one can do the same argument.
We will use a slightly weaker form of Lemma \ref{clclt}; more precisely,
we will assume that 
\begin{equation}  \label{clclt2}
\Prob \{ S _m = \frac z2 + w | S _n = z \} = 
2 \sqrt{ \frac {1}{2 \pi \sigma^2_{n,z}}} 
 \; \exp \left\{ -\frac{w^2}{2 \sigma^2_{n,z}} + 
O( \frac 1{\sqrt n}  + \frac {|w|^3}{n^2})\right\},
\end{equation}
where $\sigma^2_{n,z} = (n/4) \, [1 - (z/n)^2].$
We have replaced the $O(1/n),O(w^4/n^3)$ terms with the larger $
O(1/\sqrt n), O(|w|^3/n^2)$ terms, respectively.
Our reason for doing this is that more general random walks satisfy
local central limit theorems with this weaker error term, and it is 
useful to know that the arguments in this section only use
the weaker form.

Let $X$ denote a $N(0,1)$ random variable, let
\[  Z = Z_{n,z} = \frac{\sqrt n}{2}\,  X + \frac{z}{2},
\;\;\;\;\;  \hat Z = \hat Z_{n,z} = \sigma_{n,z}
  \, X + \frac{z}{2} , \]
and let $W= W_{n,z}$ be the random variable with distribution
of  $S_{n/2}^{(n,z)}$
that is quantile-coupled with $X$.  Note that $W$ is
also quantile-coupled with $Z$ and $\hat  Z$.  
 Let us write $F = F_{n,z}$ for the distribution function of
  $\hat Z$  and $G = G_{n,z}$ for the
distribution function of $W$.
It follows from Lemma~\ref{easy_lemma} below that there exist
  $c, \epsilon$ and $N$ 
such that for all $n > N$, for all $z \in L_n$ with $|z|/n < \epsilon$,
and all $w$ with $|w - (z/2)|/n < \epsilon$, 
\begin{equation}  \label{dist.func}
F(w - c( 1 + \frac{ (w - (z/2))^2}n )) \le G(w-) 
\le G(w) \le F(w + c( 1 + \frac{ (w - (z/2))^2}n )).
\end{equation}
It follows from (\ref{may29.1}) and (\ref{dist.func})
that the quantile-coupling satisfies
\[    |\hat Z -W| \leq c \, \left[1 + \frac{(W - (z/2))^2}{n}\right] \] 
for all $n > N$, provided that $|z|, |W - (z/2)| < \epsilon\, n$.  Also
 \[
|Z - \hat  Z|  =  \left[\, 1-\sqrt{1 - (\frac zn)^2}\, \right] \, |Z - \frac z2| 
  \leq c \, \frac{z^2}{n} , \]
provided that $|z|, |W - (z/2)| < \epsilon n$ (which implies the estimate
$|Z| \leq c' \, n$).  
 Therefore, for $n$ large enough the quantile-coupling satisfies
\[  |Z - W| \leq c \, \left[1 +  \frac{(W - (z/2))^2}{n}
   + \frac{z^2}{n} \right], \qquad \mbox{for} |z|, |W-(z/2)| \leq \epsilon \, n.\]
For any $\epsilon' > 0$, 
straightforward exponential estimates
show that there exists an $a > 0$
such that  (\ref{stephen.5}) 
holds for $|w| \geq \epsilon' \, n$. Hence
to get
   Lemma \ref{stephen.1}, it suffices to prove the following estimate.
We have written the estimate for random variables without the odd/even
parity issues of simple random walk and hence have dropped a factor
of $2$.  If we have a random variable
 that satisfies (\ref{clclt2}), i.e., that is supported
only on even or only on odd integers, we can convert it to a random variable
on all the integers by dividing the mass at $k$ equally between $k$
and $k-1$.

\begin{lemma}  \label{easy_lemma} For every $\tilde c,\tilde
\epsilon$ there exist $c_1,\epsilon_1, N_1$ such
that the following holds
for every positive integer $n > N_1$ and every $\sigma^2 \in [1/
\tilde c,\tilde c]$.
 Suppose  $S$ is an integer random variable  
 such that for every integer $|j| \leq \tilde \epsilon\,  n$,
\begin{equation}  \label{may30.2}
  \Prob\{S = j\} = \frac 1 {\sqrt{2 \pi \sigma^2n} } \, \exp\left\{
              -\frac{j^2}{2 \sigma^2n} + \delta(j) \right\}, \
\end{equation}
where
\[    |\delta(j)| \leq \tilde c \, \left[\frac{1}{\sqrt n} + \frac{|j|^3}{n^2} \right]. \]
Assume additionally that for every $a$ there exist positive 
$c,b$ such that
$\Prob\{|S| \ge a\, n\} \leq c\, e^{-
  bn}$, where $c$ and $b$ do not depend on $S$. 

Then if $G$ denotes the distribution function of $S$ and $F$
denotes the distribution function of an $N(0,\sigma^2n)$
random variable,
\begin{equation}  \label{may30.1}
         F\left(x - c_1 \,[1+ \frac{x^2}{n}]\right) \leq
         G(x-1) \leq G(x+1)\leq F\left(x + c_1 \, [1+\frac{x^2}{n}]\right), \;\;\;\;\;
          |x| \leq \epsilon_1 n . 
\end{equation}
\end{lemma}

\begin{proof}   This is a straightforward estimate of sums and integrals.
It suffices to prove (\ref{may30.1}) for integer $x$, and
by symmetry we can assume $x \geq 0$. Since we only need to 
establish the result for $n$ large when we write inequalities 
in this proof we will only be asserting that they are valid for $n$ 
large enough.

  Let $\bar F = 1-F, \bar G =
1-G, \bar \Phi = 1-\Phi$, where $\Phi$ is the distribution function of 
a $N(0,1)$.
Let
$q(y) = q(y;n,\sigma^2)
 =(2 \pi \sigma^2 n)^{-1/2} \, e^{-y^2/(2 \sigma^2 n)}$.  

We will deal with the somewhat easier 
case $|x| \leq \sqrt{3\tilde{c}}\sqrt n$, at once.
Note that \eqref{may30.2} implies $
c'/\sqrt n \leq G(x) - G(x-1) \leq c'/\sqrt n$ for
$|x| \leq \sqrt{3\tilde{c}}\sqrt n$.
By definition,
 \[  \bar G(x) =
  \sum_{j > x} \Prob\{S_n = j\} = \sum_{j > x} q(j)
    + \sum_{j > x} [\Prob\{S_n = j\} - q(j)]. \]
Also,
\[   \sum_{j > x} q(j)  =
 O(\frac{1}{\sqrt n}) +
  \int_{x}^\infty q(y) \; dy  .\]
Straightforward estimates using the assumptions give
\[  \sum_{j > x} \left|\Prob\{S_n = j\} - q(j)\right|
  = O(\frac 1 {\sqrt n} ). \]
Hence, we see that (\ref{may30.1}) holds for $|x| \leq \sqrt{3\tilde{c}}
\sqrt n$.

It remains to prove \eqref{may30.1} for $x > \sqrt{3\tilde{c}}
\sqrt n$. To make our strategy more intuitive for the reader we 
state now the simple fact about the standard normal distribution that 
lies at the heart of our lemma. The interested reader is refered to 
\cite{Mason} for many more details, and further reading on the KMT 
approximation.
For any $A >0$, there exist a $c$ and 
an $\epsilon > 0$ such that for $\sigma \sqrt{n} \le z \le \epsilon n$,
\begin{align} \label{heart}
e^{ A \frac{z^3}{n^2} }
    \bar{\Phi}(\frac z{\sigma{\sqrt n}})
\le \bar{\Phi}(\frac z{\sigma{\sqrt n}} - c \frac{z^2}{\sigma n ^{3/2}}), \\
e^{- A \frac{z^3}{n^2} }
    \bar{\Phi}(\frac z{\sigma{\sqrt n}})
\ge \bar{\Phi}(\frac z{\sigma{\sqrt n}} + c \frac{z^2}{\sigma n ^{3/2}}).\nonumber
\end{align}
 
Assume now that $ \sqrt{3\tilde{c}}\sqrt n \leq x \leq n^{5/8}$. The choice of
$5/8$ as the exponent is somewhat arbitrary. One could take any exponent 
strictly less than $2/3$ as the argument below shows.
We have
\[  \sum_{j > x} |\Prob\{S = j \} - q(j)| \leq
    O(\exp\{-c_3 n^{1/3}\}) + c \sum_{j > x} \frac{j^3}{n^2} 
  \, q(j) . \] 
Using that $x$ is an integer and that $y^3q(y)$ is decreasing for 
$y > \sqrt{3\tilde{c}}\sqrt{n}$, we see that
\[  \sum_{j > x} \frac{j^3}{n^2}  q(j) \leq 
    \int_{x}^\infty \frac{y^3}{n^2} \; \frac{1}{\sqrt{2 \pi \sigma^2 n}}
   \; e^{-y^2/(2  \sigma^2n)} \; dy  \le c \frac{1}{\sqrt n}
   \int_{\frac{x}{\sigma \sqrt n}}^\infty   z^3 \;
 e^{-z^2/2} \; dz.\]
But for $s \geq 1$,
\[  \int_s^\infty  z^3 \;
 e^{-z^2/2} \; dz
   \leq c \, s^3 \int_s^\infty e^{-z^2/2} \; dz
   \leq c \, s^2 \, e^{-s^2/2} . \]
Therefore,
\[  \sum_{j > x} |\Prob\{S = j \} - q(j)|   \leq c \; \frac{x^2}{n^{3/2}}
  \; \exp\{-\frac{x^2}{2 \sigma^2 n}\} \le 
c \frac {x^3}{n^2} \bar{\Phi}(\frac x {\sigma \sqrt n}). \]
Also, using simple estimates we obtain
\[
\sum_{j > x} q(j) = \bar{\Phi}(\frac x {\sigma \sqrt n}) (1 + O(\frac x n))
.\]
Hence, 
\[
\bar{G}(x) = \bar{\Phi}(\frac x {\sigma \sqrt n}) \exp \{O(\frac {x^3}{n^2}) \}
.\]
The result for $x$ in this range now follows from \eqref{heart}.

Assume now that $n^{5/8} \le x$.
Note that there is a constant $\bar{c}$, depending only on $\tilde{c}$,
 such that $q(y) e ^{\frac{2\tilde{c}y^3}{n^2}}$ is decreasing for 
$x \le y \le 2\bar{c}n$, and such that
\[
\int_{x-1} ^{\bar{c}n} q(y) e ^{\frac{2\tilde{c}y^3}{n^2}} dy
\le \int_{x-2} ^{2x} q(y) e ^{\frac{2\tilde{c}y^3}{n^2}} dy
.\]
From this we see that there is an $\epsilon_2$ such that if $x \le 
\epsilon_2 n $, then 
\begin{align*}
\bar{G} (x) &\le \int_x ^{\bar{c}n} q(y) e ^{\frac{2\tilde{c}y^3}{n^2}} dy
+ \bar{G} (\bar{c}n-1)
\le \int_{x-1} ^{\bar{c}n} q(y) e ^{\frac{2\tilde{c}y^3}{n^2}} dy 
\le \int_{x-2} ^{2x} q(y) e ^{\frac{2\tilde{c}y^3}{n^2}} dy.
\end{align*}
Hence, 
\[
\bar{G} (x) \le e ^{\frac{16\tilde{c}x^3}{n^2}} \bar{\Phi}
(\frac{x-2}{\sigma \sqrt{n}}).
\]
Similar arguments can be used to obtain
\[
\bar{G} (x) \ge ^{-\frac{16\tilde{c}x^3}{n^2}} \bar{\Phi}
(\frac{x+2}{\sigma \sqrt{n}}).
\]
The result for $x$ in this range follows from \eqref{heart}, and this 
concludes the proof of the lemma.

\end{proof}
 
\section{Acknowledgements}
The authors would like to thank Harrison Zhou for conversations and 
useful references on the KMT approximation.

\end{document}